\documentclass[12pt]{article}  
\usepackage{amsmath, amsfonts, amsthm, amssymb, amscd}
\usepackage[mathcal]{euscript} 
\usepackage{dsfont, pxfonts, mathrsfs, accents, verbatim}
\theoremstyle{plain}
        \newtheorem{theorem}{Theorem}[section]

\numberwithin{equation}{section}

\def\eqref#1{{(\ref{#1})}} 
\newcommand \be         {\begin{equation}}
\newcommand \ee         {\end{equation}}
\newcommand \RR     {\mathbb{R}}

\newcommand \Lcal       {\mathcal L}

\newcommand \del        {{\partial}}
\newcommand \MM     {{\mathbf M^n}}
\newcommand \Lorentz   {{\mathbf M^{n+1}}}
\newcommand \Lorentzplus   {{\mathbf M_+^{n+1}}}

\newcommand \Hcal   {\mathcal H}
\newcommand \dive   {\mbox{div}_g\hskip.06cm}

\newcommand \grade  {\mbox{grad}_g\hskip.06cm}
\newcommand \ubar   {{\bar u}}

\newcommand \lam {\lambda}

\newcommand \vol     {\mbox{vol}}
\newcommand \Riem   	{\mbox{Rm}}
\newcommand \gh     	{{\widehat g}\,}

\newcommand \inj    	{\mbox{Inj}}

\begin{document}

\title{Hyperbolic conservation laws 
\\
and spacetimes with limited regularity}

\author{\Large Philippe G. LeFloch
\\
\\ 
\sl Laboratoire Jacques-Louis Lions 
\\
\sl Centre National de la Recherche Scientifique
\\
\sl Universit\'e de Paris 6, 4, Place Jussieu, 75252 Paris, France. 
\\
\sl Email: \texttt{LeFloch@ann.jussieu.fr}
}
\date{October 29, 2006} 
\maketitle

\section{Introduction} 
 
Hyperbolic conservation laws posed on manifolds arise in many applications to 
geophysical flows and general relativity. Recent work by the author and his collaborators 
attempts to set the foundations for a study of weak solutions defined on Riemannian or Lorentzian 
manifolds and includes an investigation of the existence and qualitative behavior of solutions. 
The metric on the manifold may either be fixed (shallow water equations on the sphere, for instance)
or be one of the unknowns of the theory (Einstein-Euler equations of general relativity). This work is especially 
concerned with solutions and manifolds with limited regularity. 
We review here results on three themes: 
(1) Shock wave theory for hyperbolic conservation laws on manifolds, 
developed jointly with M. Ben-Artzi (Jerusalem); 
(2) Existence of matter Gowdy-type spacetimes with bounded variation,  
developed jointly with J. Stewart (Cambridge). 
(3) Injectivity radius estimates for Lorentzian manifolds under curvature bounds,
developed jointly with B.-L. Chen (Guang-Zhou).

\section{Conservation laws on a Riemannian manifold}

In the present section, $(\MM,g)$ is a compact, oriented, $n$-dimensional Riemannian manifold.
As usual, the tangent space at a point $x \in \MM$ is denoted by
$T_x \MM$ and the tangent bundle by $T\MM : = \bigcup_{x \in \MM} T_x \MM$,
while the cotangent bundle is denoted by $T^\star \MM = T_x^\star \MM$. 
The metric structure is determined by a positive-definite, $2$-covariant tensor field $g$.

A {\sl flux} on the manifold $\MM$ is a vector field $f=f_x(\ubar)$ depending smoothly upon the parameter $\ubar$. 
The conservation law associated with $f$ reads 
\be
\del_t u + \dive (f(u))  = 0,
\label{2.3}
\ee
where the unknown $u=u(t,x)$ is defined for $t \geq 0$ and $x \in \MM$ and the
divergence operator is applied to the vector field $x \mapsto f_x(u(t,x)) \in T_x\MM$. We say that 
the flux is {\sl geometry-compatible} if
\be
\label{58}
\dive f_x(\ubar) = 0, \qquad \ubar \in \RR, \, x \in \MM.
\ee 
We propose to single out this class of conservation laws as an important case of interest, 
which leads to robust $L^p$ estimates that do not depend on the geometry of the manifold. The equations 
arising in continuum physics do satisfy this condition. 

Equation \eqref{2.3} is a geometric partial differential equation which depends on the geometry
of the manifold only.
All estimates must take a coordinate-independent form; 
in the proofs however, it is often convenient to introduce particular coordinate charts.
We are interested in solutions ${u \in L^\infty(\RR_+\times\MM)}$ 
assuming a prescribed initial condition $u_0 \in L^\infty(\MM)$:
\be
\label{2.init} 
u(0,x) = u_0(x), \quad x \in \MM.
\ee
We extend Kruzkov theory of the
Euclidian space $\RR^n$ (see \cite{Kruzkov}) to the case of a Riemannian manifold, as follows. 

Let $f=f_x(\ubar)$ be a geometry-compatible flux on the Riemannian manifold $(M,g)$.
A {\sl convex entropy/entropy-flux pair} is a pair $(U,F)$ where $U: \RR \to \RR$
is convex and $F=F_x(\ubar)$ is the vector field defined by 
$$
F_x(\ubar) := \int^{\ubar}
\hskip-.15cm \del_{u'} U(u') \, \del_{u'} f_x(u') \, du',  \qquad \ubar \in \RR, \, x \in \MM. 
$$
Given $u_0 \in L^\infty(\MM)$ a function $u \in L^\infty\bigl(\RR_+, L^\infty(\MM)\bigr)$ is
called an {\sl entropy solution} to the initial value problem \eqref{2.3}, \eqref{2.init} if
the following entropy inequalities hold
\begin{eqnarray*} 
& \displaystyle \int\int_{\RR_+ \times \MM}  \hskip-.15cm \Bigl( U(u(t,x)) \, \del_t \theta(t,x)
+ g_x\big(F_x(u(t,x)), \grade \theta(t,x) \big) \Bigr) \, dV_g(x) dt
\\
& \displaystyle + \int_\MM  \hskip-.15cm U(u_0(x)) \, \theta(0,x) \, dV_g(x) \geq 0,
\end{eqnarray*} 
for every convex entropy/entropy flux pair $(U,F)$
and all smooth function $\theta = \theta(t,x) \geq 0$ compactly supported in $[0,\infty)\times \MM$.

\begin{theorem}[Well-posedness theory on a Riemannian manifold. I]
\label{3-1}
Let $f=f_x(\ubar)$ be a geometry-compatible flux on a Riemannian manifold $(\MM,g)$.
Given $u_0 \in L^\infty(\MM)$ there exists a unique entropy solution
$u \in L^\infty(\RR_+ \times \MM)$ to the problem \eqref{2.3}--\eqref{2.init}. Moreover, for each $1 \leq p \leq \infty$, 
$$
\|u(t) \|_{L^p(\MM; dV_g)} \leq \|u_0\|_{L^p(\MM;dV_g)}, \qquad t \in \RR_+,
$$
and, given two entropy solutions $u,v$ associated with initial data $u_0, v_0$, 
$$ 
\| v(t) - u(t) \|_{L^1(\MM; dV_g)} \leq \| v_0 - u_0 \|_{L^1(\MM; dV_g)}, \qquad t \in \RR_+.
$$
\end{theorem}

The framework proposed here allows us to construct entropy solutions on a Riemannian manifold 
via the vanishing diffusion method or the finite volume method  \cite{AmorimBenArtziLeFloch,BenArtziLeFloch}.  
Following DiPerna~\cite{DiPerna} we can introduce the (larger) class of entropy measure-valued solutions
$(t,x) \in \RR_+ \times \MM \mapsto \nu_{t,x}$.

\begin{theorem}[Well-posedness theory on a Riemannian manifold. II] 
\label{3-4}
Let $f=f_x(\ubar)$ be a geometry-compatible flux on a Riemannian manifold $(\MM,g)$.
Let $\nu$ be an entropy measure-valued solution
to \eqref{2.3}--\eqref{2.init} for some $u_0 \in L^\infty(\MM)$.
Then, for almost every $(t,x)$, $\nu_{t,x} = \delta_{u(t,x)}$, 
where $u \in L^\infty(\RR_+ \times \MM)$ is the unique entropy solution to the problem. 
\end{theorem}

Finally we can relax the geometry compatibility condition and 
consider a general conservation law associated with an arbitrary flux $f$.
More general conservation laws solely enjoy the $L^1$ contraction property and leads to a unique contractive
semi-group of entropy solutions.  

\begin{theorem}[Well-posedness theory on a Riemannian manifold. III]
Let $f=f_x(\ubar)$ be an arbitrary (not necessarily divergence-free) flux on $(\MM,g)$,
satisfying the linear growth condition
$$
| f_x(\ubar) |_g \lesssim 1 + |\ubar|, \qquad \ubar \in \RR, \, x \in \MM.
$$
Then there exists a unique contractive, semi-group of entropy solutions
$u_0 \in L^1(\MM) \mapsto u(t) := S_t u_0 \in L^1(\MM)$
to the initial value problem \eqref{2.3}, \eqref{2.init}.
\end{theorem}
 
For the proofs we refer to \cite{BenArtziLeFloch}. See \cite{AmorimBenArtziLeFloch} for 
the convergence of the finite volume schemes on a manifold. 
Earlier material can be found in Panov \cite{Panov} ($n$-dimensional manifold) 
and in LeFloch and Nedelec \cite{LeFlochNedelec} (Lax formula for general metrics including 
the case of spherical symmetry).


\section{Conservation laws on a Lorentzian manifold}

Motivated by the application to general relativity, we can extend the theory to a Lorentzian manifold.
Let $(\Lorentz,g)$ be a time-oriented, $(n+1)$-dimensional Lorentzian manifold,
$g$ being a metric tensor with signature $(-, +, \ldots, +)$.
Tangent vectors $X$ can be separated into time-like vectors
($g(X,X) < 0$), null vectors ($g(X,X) = 0$), and space-like vectors ($g(X,X) > 0$). 
The null cone separates time-like vectors into future-oriented and past-oriented ones. 
Let $\nabla$ be the Levi-Cevita connection associated with the Lorentzian metric $g$.

A {\sl flux} on the manifold $\Lorentz$
is a vector field $x \mapsto f_x(\ubar) \in T_x\Lorentz$, depending on a parameter $\ubar \in \RR$.
The conservation law on $(\Lorentz,g$) associated with $f$ is
\be
\dive \big( f(u) \big) = 0, \qquad u: \Lorentz \to \RR.
\label{LO.1b}
\ee
It is said to be {\sl geometry compatible} if 
\be
\dive f_x (\ubar) = 0, \quad \ubar \in \RR, \, x \in \Lorentz.
\label{LO.divergencefree}
\ee
Furthermore, $f$ is said to be a {\sl time-like flux} if
$g_x\big(\del_u f_x(\ubar), \del_u f_x(\ubar) \big) < 0$, $x \in \Lorentz$, $\ubar \in \RR$.

Note that our terminology here differs from the one in the Riemannian case,
where the conservative variable was singled out.
We are interested in the initial-value problem associated with \eqref{LO.1b}. We fix a space-like hypersurface
$\Hcal_0 \subset \Lorentz$ and a measurable and bounded function $u_0$ defined on $\Hcal_0$.
Then, we search for $u=u(x) \in L^\infty(\Lorentz)$ satisfying \eqref{LO.1b}  in the distributional
sense and such that the (weak) trace of $u$ on $\Hcal_0$ coincides with $u_0$:
\be
u_{|\Hcal_0} = u_0.
\label{LO.initial}
\ee
It is natural to require that the vectors $\del_u f_x(\ubar)$ are time-like and future-oriented.

We assume that the manifold $\Lorentz$ is {\sl globally hyperbolic,} in the sense that
there exists a foliation of $\Lorentz$ by space-like, compact, oriented hypersurfaces $\Hcal_t$ ($t \in \RR$):
$\Lorentz = \bigcup_{t \in \RR} \Hcal_t$. 
Any hypersurface $\Hcal_{t_0}$ is referred to as a Cauchy surface in $\Lorentz$,
while the family $\Hcal_t$ ($t \in \RR$) is called an admissible foliation associated with $\Hcal_{t_0}$.
The future of the given hypersurface will be denoted by
$\Lorentzplus := \bigcup_{t \geq 0} \Hcal_t$. 
Finally we denote by $n^t$ the future-oriented, normal vector field to each $\Hcal_t$,
and by $g^t$ the induced metric.
Finally, along $\Hcal_t$, we denote by $X^t$ the normal component of a vector field $X$, thus
$X^t := g(X,n^t)$.
 
A flux $F=F_x(\ubar)$ is called a {\sl convex entropy flux} associated with the conservation law \eqref{LO.1b}
if there exists a convex $U:\RR \to \RR$ such that
$$
F_x(\ubar) = \int_0^\ubar \del_u U(u') \,\del_u f_x(u') \, du', \qquad x \in \Lorentz, \, \ubar \in \RR.
$$
A measurable and bounded function $u=u(x)$ is called an {\sl entropy solution}
of the geometry-compatible conservation law \eqref{LO.1b}-\eqref{LO.divergencefree} if
$$
\int_{\Lorentzplus} g(F(u), \grade_g \theta) \, dV_g + \int_{\Hcal_0} g_0(F(u_0), n_0) \, \theta_{\Hcal_0} \, dV_{g_0} \geq 0.
$$
for all convex entropy flux $F=F_x(\ubar)$ and all smooth $\theta \geq 0$ compactly supported in
$\Lorentzplus$.

\begin{theorem} [Well-posedness theory on a Lorentzian manifold]
\label{LO-theo}
Consider a geometry-compatible conservation law \eqref{LO.1b}-\eqref{LO.divergencefree}
posed on a globally hyperbolic Lorentzian manifold $\Lorentz$.
Let $\Hcal_0$ be a Cauchy surface in $\Lorentz$, and $u_0: \Hcal_0 \to \RR$ be measurable and
bounded. Then, the initial-value problem \eqref{LO.1b}-\eqref{LO.initial}
admits a unique entropy solution $u=u(x) \in L^\infty(\Lorentz)$.
For every admissible foliation $\Hcal_t$, the trace $u_{\Hcal_t}$ exists and belong to $L^1(\Hcal_t)$,
and 
$\| F^t(u_{\Hcal_t} \|_{L^1(\Hcal_t)}$
is non-increasing in time, for any convex entropy flux $F$. Moreover, given any two entropy solutions $u,v$,
$$
\| f^t(u_{\Hcal_t}) - f^t(v_{|\Hcal_t}) \|_{L^1(\Hcal_t)}
$$
is non-increasing in time.
\end{theorem}

We emphasize that, in the Lorentzian case, no time-translation property is available in general, contrary to the
Riemannian case. Hence, no time-regularity is implied by the $L^1$ contraction property.


\section{Existence \ of \ matter \ Gowdy-type \ spacetimes \ with \ bounded \ variation}

Vacuum Gowdy spacetimes are inhomogeneous spacetimes admitting
two commuting spatial Killing vector fields. The existence of 
{\sl vacuum} spacetimes with Gowdy symmetry is well-known and 
the long-time asymptotics of solutions have been found to be particularly complex.  
In comparison, much less emphasis has been put on {\sl matter} spacetimes.   
Recently, LeFloch, Stewart and collaborators \cite{BLSS,LeFlochStewart} initiated a rigorous mathematical 
treatment of the coupled Einstein-Euler system on Gowdy spacetimes.  
The unknowns of the theory are the density and velocity of the fluid together 
with the components of the metric tensor. The existence for the Cauchy problem 
in the class of solutions with (arbitrary large) bounded total variation is proven by a generalization 
of the Glimm scheme. Our theory allows for the formation of shock waves in the fluid and singularities in the geometry.  
The first results on shock waves and the Glimm scheme 
in special and general relativity are due to Smoller and Temple \cite{SmollerTemple} 
(flat Minkowski spacetime) and Groah and Temple \cite{GroahTemple}
(spherically symmetric spacetimes).  
The novelty in \cite{BLSS,LeFlochStewart} is the generalization to a model 
allowing for both gravitational waves and shock waves.  

The metric is given in the polarized Gowdy symmetric form   
\be
ds^2 = e^{2a} \, (- dt^2 + dx^2) 
       + e^{2b} \, (e^{2c} \, dy^2 + e^{-2c} \, dz^2), 
\label{2.1}
\ee 
where the variables $a,b,c$ depend on the time variable $t$ and the space variable $x$, only.  
We consider {\sl Einstein field equations} $G^{\alpha\beta} = \kappa T^{\alpha\beta}$ 
for perfect fluids with energy density $\mu>0$ and pressure  
$p = \mu {c_s^2}$. Here, the sound speed $c_s$ is a constant with $0<c_s<1$  
and $G^{\alpha\beta}$ denotes the Einstein tensor, while $\kappa$ is a normalization constant. 

The $4$-velocity vector $u^\alpha$ of the fluid 
is time-like and is normalized to be of unit length
and we define the scalar velocity $v$ and relativistic factor $\xi=\xi(v)$ by   
$(u^\alpha) = e^{-a} \, \xi \, (1,v,0,0)$
and $\xi = (1-v^2)^{-1/2}$. 
The matter is described by the energy-momentum tensor
$T^{\alpha\beta} = (\mu + p) \, u^\alpha u^\beta + p \, g^{\alpha\beta}$, 
from which we extract the fields $\tau$, $S$ and $\Sigma$: 
\begin{eqnarray*}
T^{00} = & e^{-2a} \big( (\mu+p) \, \xi^2 - p \bigr) & =: e^{-2a}\tau, 
\\
T^{01} = & e^{-2a}  (\mu+p) \, \xi^2 v & =: e^{-2a} S, 
\\
T^{11} = & e^{-2a} \big( (\mu+p) \xi^2 v^2 + p \bigr) & =: e^{-2a} \, \Sigma.
\end{eqnarray*} 

After very tedious calculations we arrive at the \textbf{constraint equations} 
\begin{eqnarray*}
& 2 \, a_t \, b_t + 2 \, a_x \, b_x + b_t^2 - 2 \, b_{xx} - 3 \, b_x^2
- c_t^2 - c_x^2 & = \kappa \, e^{2a} \, \tau, 
\\
& -2 \, a_t \, b_x - 2 \, a_x \, b_t + 2 \, b_{tx} + 2 \, b_t \, b_x + 
2 \, c_t \, c_x  & = \kappa \, e^{2a} \, S, 
\end{eqnarray*}
and the \textbf{evolution equations}  
\begin{eqnarray*}
& a_{tt} - a_{xx} & = b_t^2 - b_x^2 - c_t^2 + c_x^2 +
     \frac{\kappa}{2} \, e^{2a} \, ( - \tau + \Sigma -2 \, p ),
\\
& b_{tt} - b_{xx} & = -2 b_t^2 + 2 b_x^2 + 
 \frac{\kappa}{2} \, e^{2a} \, ( \tau - \Sigma),
\\
& c_{tt} - c_{xx} & = -2 \, b_t \, c_t + 2 \, b_x \, c_x.
\end{eqnarray*}

The evolution equations for the fluid, $\nabla_\beta T^{\alpha\beta}=0$, are the \textbf{Euler equations} 
$$
\tau_t + S_x = T_1, \qquad 
S_t + \Sigma_x = T_2,
$$
in which the source terms $T_1, T_2$ are nonlinear in first-order derivatives of the metric
and fluid variables.

We propose to reformulate the Einstein-Euler equations in the form of a nonlinear
hyperbolic system of balance laws with integral source-term, in the variables $(\mu, v)$ and 
$w := \bigl(a_t, a_x, \beta_t, \beta_x, c_t, c_x\bigr)$, 
where  $\beta = e^{2b}$. It is convenient to also set $\alpha = e^{2a}$.  
The functions $\alpha, b$ (and $a, \beta$) are determined by 
\begin{eqnarray}
& \alpha(t,x) = e^{2 a(t,x)},
\qquad a(t,x) = \int_{-\infty}^x w_2(t,y) \, dy, 
\nonumber 
\\
& b (t,x) = \frac{1}{2} \, \ln \beta(t,x), \qquad 
\beta(t,x) = 1+ \int_{-\infty}^x w_4(t,y) \, dy. 
\nonumber
\end{eqnarray}
Obviously, we are interested in solutions such that $\beta$ remains positive. 

The equations under consideration consist of three sets of two equations  
associated with the propagation speeds $\pm 1$, the speed of light (after normalization).  
The principal part of the fluid equations are the standard relativistic fluid
equations in  a Minkowski background, with wave speeds 
$\lam_\pm = (v \pm c_s) / (1 \pm v \, c_s)$. 
To formulate the initial-value problem it is natural to prescribe the values 
of $\mu,v,w$ on the initial hypersurface at $t=0$, denoted by $(\mu^0,v^0,w^0)$. 

Our main result is: 

\begin{theorem}[Existence of Gowdy spacetimes with compressible matter]
 \label{3-1b}
Consider the $(\mu,v,w)$-formulation of the Einstein-Euler equations on a polarized 
Gowdy spacetime with plane-symmetry.
Let the initial data $(\mu^0,v^0,w^0)$ be of bounded total variation,  
$TV(\mu^0,v^0,w^0) < \infty$, satisfying the constraints, and suppose that the corresponding 
functions $\alpha^0,b^0$ are measurable and bounded,  $\sup |(\alpha^0,b^0)| < \infty$.  
Then the Cauchy problem admits a weak solution $\mu,v,w$ such that for some increasing $C(t)$  
$$
TV(\mu,v,w)(t) + \sup|(\alpha,b)(t,\cdot)| 
\leqslant C(t), \qquad t \geq 0,  
$$  
and are defined up to a maximal time $T \leqslant \infty$. 
If $T < \infty$ then either the geometry variables $\alpha,b$  blow up: 
$\lim_{t \to T} \bigl( \sup_\RR |\alpha(t,\cdot)| + |b(t,\cdot)| \bigr) = \infty$, 
or the energy density blows up :  
$\lim_{t \to T} \sup_\RR |\mu(t,\cdot)| = \infty$.  
\end{theorem}

Hence, the solution exists until either a singularity occurs in the geometry
(e.g.~the area $\beta$ of the $2$-dimensional space-like orbits of the
symmetry group vanishes) or the matter collapses to a point. 
To our knowledge this is the first {\sl global} existence result for the Euler-Einstein equations. 

If a shock wave forms in the fluid, then $\mu, v$ will be discontinuous and, as a consequence,  
$w_{3x}$ and $w_{4x}$ might also be discontinuous.  
In fact, Theorem~5 allows not only such discontinuities in second-order
derivatives of the geometry components (i.e.~at the level of the curvature), 
but also discontinuities in the first-order derivatives which propagate at the speed of light.  
The latter correspond to {\sl Dirac distributions} in the curvature of the metric.


\section{Lower bounds on the injectivity radius of Loren\-tzian manifolds}

Motivated by the application to spacetimes of general relativity 
and by earlier results established by Anderson \cite{Anderson}
and Klainerman and Rodnianski \cite{KR}, 
we investigate in \cite{ChenLeFloch} 
the geometry and regularity of $(n+1)$-dimensional Lorentzian manifolds $(M,g)$. 
Under curvature and volume bounds we establish new injectivity radius estimates which 
are valid either in arbitrary directions or in null cones. Our estimates are purely local and are formulated 
via the ``reference'' Riemannian metric $\gh_T$ associated with an arbitrary future-oriented time-like vector field
$T$. 

Our proofs are based on suitable generalizations of arguments from Riemannian geometry 
and rely on the observation that geodesics in the Euclidian and Minkowski spaces
coincide, so that estimates for the reference Riemannian metric can be carried over to 
the Lorentzian metric. Our estimates should be useful to investigate the qualitative behavior of spacetimes satisfying Einstein field equations. 
    
We state here one typical result from \cite{ChenLeFloch} encompassing a large class of Lorentzian manifolds. 
Fix a point $p\in M$, and let us assume that a domain $\Omega \subset M$ containing $p$ 
is foliated by spacelike hypersurfaces $\Sigma_t$ with normal $T$,
say $\Omega = \bigcup_{t \in [-1,1]} \Sigma_t$. 
Assume also that the geodesic ball $B_{\Sigma_0}(p,1) \subset \Sigma_0$  
is compactly contained in $\Sigma_0$. Consider the following assumptions
where $K_0, K_1, K_2$ and $v_0$ are positive constants: 
\be
\label{A1}
K_0 \leq - \big| \frac{\del}{\del t} \Big|_g^2 \leq 1/K_0 \quad \mbox{ in } \Omega,
\ee
\be
\label{A2}
|\Lcal_Tg |_{\gh_T} \leq K_1 \quad \mbox{ in } \Omega,
\ee
\be
\label{A3} 
|\Riem_g |_{\gh_T} \leq K_2 \quad  \mbox{ in } \Omega, 
\ee
\be
\label{A4} 
\vol_g(B_{\Sigma_0}(p,1))\geq v_0,
\ee  

We prove in \cite{ChenLeFloch} : 

\begin{theorem}[Injectivity radius estimate for Lorentzian manifolds] 
\label{inject}
Let $(M,g)$ be a Lorentzian manifold satisfying \eqref{A1}--\eqref{A4} at a point $p \in M$. Then, there exists 
a positive constant $i_0$ depending on the foliation bounds $K_0,K_1$, curvature bound $K_2$, 
volume bound $v_0$, and dimension $n$ so that the injectivity radius at $p$ is bounded below by $i_0$, that is 
$\inj(M,g,p) \geq i_0$. 
\end{theorem}

 
\section*{Acknowledgments} 
The author was partially supported 
by the A.N.R. grant 06-2-134423 entitled ``Mathematical Methods in General Relativity'' (MATH-GR).

 \end{document}